\documentclass{article}

\usepackage{amsmath}
\usepackage{amsfonts}
\usepackage{amssymb}
\usepackage{amsthm}
\usepackage{url}

\renewcommand\iff{\leftrightarrow}
\newcommand\imp{\rightarrow}
\newcommand\logand{\wedge}
\newcommand\logor{\vee}

\newcommand\three{\mathbf{3}}

\newcommand\simp{\Rightarrow}
\newcommand\lps{\ensuremath{\text{LP}_{\simp}}}

\newcommand\siff{\Leftrightarrow}
\newcommand\gbot{{>\bot}}
\newcommand\con{{\text{con}}}
\newcommand\entails{\vdash}

\newcommand\true{\text{T}}

\newcommand\tr{\text{Tr}}
\newcommand\cl{\text{cl}}

\theoremstyle{definition}
\newtheorem{definition}{Definition}
\newtheorem{theorem}{Theorem}

\usepackage{natbib}
\title{\lps: Extending LP with a strong conditional operator.}
\author{Nick Thomas}

\begin{document}

\maketitle

\begin{abstract}
We augment LP with a strong conditional operator, to yield a logic we call ``strong LP,'' or \lps. The resulting logic can speak of consistency in more discriminating ways, but introduces new possibilities for trivializing paradoxes.
\end{abstract}

\section{Introduction}
\label{sec:intro}

A popular paraconsistent logic is the logic LP (``Logic of Paradox''), introduced by \citet{asenjo66} and discussed, e.g., in \citet{priest79}, \citet{priest84}, and \citet{priest02}. The idea behind LP is to take classical logic's two-valued semantics and extend it to a three-valued semantics, with truth values true ($\top$), false ($\bot$), and paradoxical ($+$). We think of the third truth value as meaning ``both true and false.'' We extend the classical truth tables as in Table~\ref{tbl:lp-operators}.

\begin{table}
\begin{center}
\begin{tabular}{ccc}

\begin{tabular}{c|c}
       & $\neg$ \\
\hline
$\top$ & $\bot$ \\
$+$    & $+$ \\
$\bot$ & $\top$
\end{tabular} &
\begin{tabular}{c|c|c|c}
$\logand$ & $\top$ & $+$    & $\bot$ \\
\hline
$\top$    & $\top$ & $+$    & $\bot$ \\
$+$       & $+$    & $+$    & $\bot$ \\
$\bot$    & $\bot$ & $\bot$ & $\bot$
\end{tabular} &
\begin{tabular}{c|c|c|c}
$\logor$ & $\top$ & $+$    & $\bot$ \\
\hline
$\top$   & $\top$ & $\top$ & $\top$ \\
$+$      & $\top$ & $+$    & $+$    \\
$\bot$   & $\top$ & $+$    & $\bot$
\end{tabular} \\ \\

&
\begin{tabular}{c|c|c|c}
$\imp$ & $\top$ & $+$    & $\bot$ \\ \hline
$\top$     & $\top$ & $+$    & $\bot$ \\
$+$        & $\top$ & $+$    & $+$    \\
$\bot$     & $\top$ & $\top$ & $\top$
\end{tabular} &
\begin{tabular}{c|c|c|c}
$\iff$   & $\top$ & $+$    & $\bot$ \\ \hline
$\top$   & $\top$ & $+$    & $\bot$ \\
$+$      & $+$    & $+$    & $+$    \\
$\bot$   & $\bot$ & $+$    & $\top$
\end{tabular}

\end{tabular}
\end{center}
\caption{Truth tables for LP's logical operators.}
\label{tbl:lp-operators}
\end{table}

Why are the truth tables defined as they are? One natural way to answer is to give an equivalent definition, to be found e.g. in \citet{priest02}. Let us define $\top = 1, + = .5, \bot = 0$. Then we may define $a \logand b = \min\{a,b\}$, $a \logor b = \max\{a,b\}$, and $\neg a = 1 - a$. This definition hopefully shows some of the symmetries present in the definition.

LP validates all classical tautologies; however, it fails to validate all classical inferences. \citep{priest02} In particular, modus ponens fails; $p \imp q, p \nvDash q$. This, arguably, forms a serious obstacle to using LP for much of anything. But LP is a very simple and intuitive proposal, and the obstacle has not made researchers give up on it. Various ways around the problem have been proposed.

\citet{priest91} proposes a system called ``minimally inconsistent LP,'' or $\text{LP}_{\text{m}}$, wherein we restrict our attention to ``minimally inconsistent'' models: essentially, models which satisfy as few contradictions as possible. Then we say that $T \models U$ iff every minimally inconsistent model of $T$ satisfies $U$. This lets us recover all techniques of classical reasoning in the consistent case. \citet{crabbe13} has recently done some important work on this theory.

\citet{beall11} proposes another solution. Beall defines a system called ``multiple conclusion LP,'' wherein we say that $T$ satisfies $U$ iff there is no interpretation which values all members of $T$ as $\top$ or $+$, and values all members of $U$ as $\bot$. The idea is that the set of conclusions represents a set of alternatives. We interpret $p, p \imp q \models q, p \logand \neg p, q \logand \neg q$ as meaning ``if $p$ and $p \imp q$, then either $q$, or one of our premises is inconsistent.'' This lets us get back classical logic by adding appropriate contradictions to our set of alternative conclusions.

A third solution, due to \citet{beall13}, is to add rules of the form $p \entails \bot$, called ``shrieks,'' where such a rule expresses that $p$ is not a theorem, up to triviality. This gives a way to say that a theorem is ``true and not false,'' or non-paradoxically true; we just assert $p$ and $\neg p \entails \bot$. Then all forms of classical reasoning become valid when we are talking about shrieked propositions.

The solution we use here is to augment LP with a new logical operator. It is an alternative form of implication, which we call ``strong implication,'' and write $\simp$. It has the following truth table:
\medskip
\begin{center}
\begin{tabular}{c|c|c|c}
$\simp$    & $\top$ & $+$    & $\bot$ \\ \hline
$\top$     & $\top$ & $\bot$ & $\bot$ \\
$+$        & $\top$ & $\top$ & $\bot$ \\
$\bot$     & $\top$ & $\top$ & $\top$
\end{tabular}
\end{center}
\medskip
In the real number interpretation of truth values, we define
\[ a \simp b = \begin{cases}
1 & a \leq b \\
0 & a > b.
\end{cases} \]
We refer to LP augmented with the strong implication operator as ``strong LP,'' or $\lps$. This operator has a number of nice properties, which we shall discuss in Section~\ref{sec:lps-properties}. It implements modus ponens. Its addition makes $\lps$ functionally complete. It obeys contraction and contraposition, allows a form of proof by contradiction, and obeys restricted forms of weakening and the deduction theorem. The biconditional $p \siff q$ is true just in case $p$ and $q$ have the same truth value.

A further nice property of $\lps$ is that we can reproduce classical logic in it. As discussed in Section~\ref{sec:lps-properties}, if an $\lps$ theory has only consistent models, then its $\lps$ consequences are precisely its classical consequences; and, if an $\lps$ theory has a consistent model, then we can add to it a nontrivializing axiom schema which makes its $\lps$ consequences become exactly its classical consequences.

Additionally, we can use $\simp$ to define a number of other useful logical operators. These operators let us talk about inconsistency in very fine and discriminating ways. We can express that a statement is non-paradoxically true; that a statement is either true or paradoxical, but not false; and so forth.

The first operators are $p^\top$ (``$p$ is true''), $p^+$ (``$p$ is paradoxical''), and $p^\bot$ (``$p$ is false''). These operators yield true if $p$ has the superscripted truth value, and false otherwise. We also define $p^\gbot$ (``$p$ is not false''), which is true unless $p$ is non-paradoxically false, and $p^\con$ (``$p$ is consistent''), which yields true unless $p$ is paradoxical. The truth tables of these operators are given in Table~\ref{tbl:defined-consistency-ops} (along with the truth table for $\simp$, for comparison). Table~\ref{tbl:consistency-op-definitions} states how they may be defined in terms of $\simp$. The possibility of defining these operators shows that $\simp$ adds a great deal of power to the language. But in a sense, it adds too much power, as we shall now see.

\begin{table}
\begin{equation}
\begin{array}{ccccc}
p^\top       & p^+       & p^\bot       & p^\gbot   & p^\con \\
\top \simp p & p \siff + & p \simp \bot & + \simp p & p \simp p^\top
\end{array}
\end{equation}
\caption{Definitions of extended unary operators.}
\label{tbl:consistency-op-definitions}
\end{table}

\begin{table}
\begin{equation}
\begin{array}{ccccc}
&
\begin{array}{c|ccc}
\simp & \top & +   & \bot \\
\hline
\top  & \top & \bot & \bot \\
+     & \top & \top & \bot \\
\bot  & \top & \top & \top
\end{array} & &

\begin{array}{c|ccc}
\siff & \top & +   & \bot \\
\hline
\top  & \top & \bot & \bot \\
+     & \bot & \top & \bot \\
\bot  & \bot & \bot & \top
\end{array} & \\ \\

\begin{array}{c|c}
p & p^\top \\
\hline
\top & \top \\
+    & \bot \\
\bot & \bot
\end{array} &

\begin{array}{c|c}
p & p^+ \\
\hline
\top & \bot \\
+    & \top \\
\bot & \bot 
\end{array} &

\begin{array}{c|c}
p & p^\bot \\
\hline
\top & \bot \\
+    & \bot \\
\bot & \top
\end{array} &

\begin{array}{c|c}
p & p^\gbot \\
\hline
\top & \top \\
+    & \top \\
\bot & \bot
\end{array} &

\begin{array}{c|c}
p & p^\con \\
\hline
\top & \top \\
+    & \bot \\
\bot & \top
\end{array}
\end{array}
\end{equation}

\caption{Truth tables for extended logical operators.}
\label{tbl:defined-consistency-ops}
\end{table}

The idea of augmenting LP with additional logical connectives is not a new one; see, e.g., \citet{denyer89} and \citet{priest89}. Researchers have largely rejected it, because when combined with tools of self reference such as the $T$-schema or the na\"ive comprehension schema, it tends to yield paradoxes which lead to triviality. For instance, define the following variation of the Russell set:

\begin{equation}
\label{form:super-russell}
R = \{x : (x \in x)^\bot\}. 
\end{equation}

The statement $R \in R$ for this ``super-Russell'' set gives us a more vicious type of paradox, which is trivializing in $\lps$, as the reader may verify. A similar, Curry-like paradox arises if we set up a $T$-schema in $\lps$.

In the set-theoretic case, the simple solution is to disallow the use of $\simp$ inside set-builders. This is undoubtedly a sacrifice, but \emph{prime facie}\/ appears to leave us with a usable set theory. In the case of the $T$-schema, such a restriction would seem to contradict the very point of the schema, and so we can only say that $\lps$ is not an appropriate setting for such a schema.

\section{The logic \lps}
\label{sec:lps}

We proceed with defining $\lps$. We assume the usual definition of signatures, with constant and relation symbols but no function symbols. We assume an infinite set $\mathcal{V}$ of variable names. Terms are variables or constants. We let $\tau, \upsilon, ...$ denote terms. Atomic formulas are logical constants ($\top, +, \bot$), relations $R(\tau_1, ..., \tau_n)$, or equalities $\tau = \upsilon$. Formulas are built up from atomic formulas using the connectives $\neg, \logand, \simp$ and quantifiers $\forall x (\phi)$. We let $\logor,\imp,\iff,\exists$ be defined in terms of $\neg, \logand, \forall$ in the usual way. We let $\phi, \psi, \zeta, ...$ denote formulas.

\begin{definition}
A ``model'' is a pair $(M,I)$ associated with a signature, where:

\begin{enumerate}
\item $M$ is a nonempty set (the universe of objects).
\item $I$ is the interpretation function, giving interpretations to constant and relation symbols. Its domain consists of all constant and relation symbols in the signature.
\item For each constant symbol $c$ in the signature, $I(c) \in M$ is an object.
\item For each $n$-ary relation symbol $R$ in the signature, $I(R) : M^n \rightarrow \three$ is an $n$-ary three-valued relation on $M$. Using currying notation, we abbreviate $I(R)(y_1,...,y_n)$ to $I(R,y_1,...,y_n)$.
\end{enumerate}

Informally, we refer to the pair $(M,I)$ as just $M$.
\end{definition}

\begin{definition}
If $(M,I)$ is a model, we say that $(M,I)$ is ``consistent'' iff, for every relation $R$ in the signature, the image of $I(R)$ does not contain $+$.
\end{definition}

\begin{definition}
A ``variable assignment'' $A$ for a model $M$ is a function $A : \mathcal{V} \rightarrow M$ which gives values to the variables.
\end{definition}

\begin{definition}
Given a model $(M,I)$ and a variable assignment $A$, we define the valuation function $\text{val}^M_A$ from terms to objects as follows. We omit the superscript and subscript where unamibiguous.

\begin{enumerate}
\item $\text{val}(c) = I(c)$, where $c$ is a constant symbol in $M$'s signature.
\item $\text{val}(x) = A(x)$, where $x \in \mathcal{V}$ is a variable.
\end{enumerate}
\end{definition}

\begin{definition}
Given a model $M$ and a variable assignment $A$, we define the truth function $\true^M_A(\phi)$ from formulas to truth values as follows. We omit the superscript and subscript where unambiguous.

\begin{enumerate}
\item $\true(R(\tau_1, ..., \tau_n)) = I(R, \text{val}(\tau_1), ..., \text{val}(\tau_n))$.
\item $
\true(\tau = \upsilon) = \begin{cases}
\top & \text{if}\ \text{val}(\tau) = \text{val}(\upsilon); \\
\bot & \text{otherwise}.
\end{cases}$
\item $\true(\neg\phi) = \neg\true(\phi)$.
\item $\true(\phi \logand \psi) = \true(\phi) \logand \true(\psi)$.
\item $\true(\phi \simp \psi) = \true(\phi) \simp \true(\psi)$.
\item $\true_A(\forall x (\phi)) = \underset{y \in M}{\bigwedge} \true_{A \lbrack x \mapsto y \rbrack}(\phi)$.\footnote{$A \lbrack x \mapsto y \rbrack$ denotes $A$ modified so that $x$ gets the value $y$.}
\end{enumerate}
\end{definition}

\begin{definition}
Given theories $T,U$, a model $M$, and a variable assignment $A$, we say:

\begin{enumerate}
\item $M,A \models T$ iff $\true^M_A(\phi) > \bot$ for all $\phi \in T$.
\item $M \models T$ iff $M,A \models T$ for all variable assignments $A$. $\models T$ iff $M \models T$ for all models $M$ in the signature of $T$.
\item $T \models U$ iff, for all models $M$, if $M,A \models T$ for all variable assignments $A$, then $M,A \models U$ for all variable assignments $A$.
\end{enumerate}

We shall write $T \models_\simp U$ when we need to be clear that we are talking about $\lps$ satisfaction, and e.g. write $\models_\cl$ to denote classical satisfaction.
\end{definition}

\section{Properties of \lps}
\label{sec:lps-properties}

Unlike classical logic, LP is functionally incomplete, in the sense that not every truth function $f : \three^n \imp \three$ may be expressed as a propositional formula of LP. We prove that $\lps$ is functionally complete, in the same sense.

\begin{table}
\begin{center}
\begin{tabular}{cccc}

\begin{tabular}{c|c}
$p$ & $p$ \\
\hline
$\top$ & $\top$ \\
$+$    & $+$    \\
$\bot$ & $\bot$
\end{tabular} &

\begin{tabular}{c|c}
$p$ & $p \logor \neg p$ \\
\hline
$\top$ & $\top$ \\
$+$    & $+$    \\
$\bot$ & $\top$
\end{tabular} &

\begin{tabular}{c|c}
$p$ & $\neg p$ \\
\hline
$\top$ & $\bot$ \\
$+$    & $+$    \\
$\bot$ & $\top$
\end{tabular} &

\begin{tabular}{c|c}
$p$ & $p \logand \neg p$ \\
\hline
$\top$ & $\bot$ \\
$+$    & $+$    \\
$\bot$ & $\bot$
\end{tabular} \\ \\

\begin{tabular}{c|c}
$p$ & $p \logand +$ \\
\hline
$\top$ & $+$ \\
$+$    & $+$    \\
$\bot$ & $\bot$
\end{tabular} &

\begin{tabular}{c|c}
$p$ & $\neg p \logor +$ \\
\hline
$\top$ & $+$ \\
$+$    & $+$    \\
$\bot$ & $\top$
\end{tabular} &

\begin{tabular}{c|c}
$p$ & $\neg p \logand +$ \\
\hline
$\top$ & $\bot$ \\
$+$    & $+$    \\
$\bot$ & $+$
\end{tabular} \\ \\

\begin{tabular}{c|c}
$p$ & $\top$ \\
\hline
$\top$ & $\top$ \\
$+$    & $\top$    \\
$\bot$ & $\top$
\end{tabular} &

\begin{tabular}{c|c}
$p$ & $+$ \\
\hline
$\top$ & $+$ \\
$+$    & $+$    \\
$\bot$ & $+$
\end{tabular} &

\begin{tabular}{c|c}
$p$ & $\bot$ \\
\hline
$\top$ & $\bot$ \\
$+$    & $\bot$    \\
$\bot$ & $\bot$
\end{tabular} &

\end{tabular}
\end{center}
\caption{Unary truth functions definable in LP.}
\label{tbl:lp-definable-unary-operators}
\end{table}

\begin{theorem}
LP (with logical constants) is functionally incomplete.

\textbf{Proof.} We claim that the set of truth functions listed in Table~\ref{tbl:lp-definable-unary-operators} is the set of unary truth functions definable in LP with logical constants. Each is definable in LP, via the listed formula.

Every atomic propositional formula defines either the identity function or a constant function; and those are in the table. I wrote a computer program to check that the negation of a truth function in the table is in the table; and to check that the conjunction of two truth functions in the table is in the table. By induction on formulas, this shows that every unary truth function definable in LP is in the table. But the table contains 10 functions, whereas there are $3^3 = 27$ unary truth functions. \qed
\end{theorem}

\begin{theorem}
$\lps$ (with logical constants) is functionally complete.

\textbf{Proof.} Consider a truth function $f :\three^n \rightarrow \three$. Given any sequence of truth values $v_1, ..., v_n \in \three$ (abbreviated $(v_i)$), let $u = f(v_1, ..., v_n)$ and define

\begin{equation}
\phi_{(v_i)}(p_1, ..., p_n) = (p_1 \siff v_1 \logand \cdots \logand p_n \siff v_n) \logand u.
\end{equation}

It is easy to see that:

\begin{equation}
\phi_{(v_i)}(p_1, ..., p_n) = \begin{cases}
u & \text{if}\ p_1, ..., p_n = v_1, ..., v_n; \\
\bot & \text{otherwise}.
\end{cases}
\end{equation}

Now define $\mathcal{V}$ as the set of sequences $v_1, ..., v_n \in \three$, and define

\begin{equation}
\psi(p_1, ..., p_n) = \underset{(v_i) \in \mathcal{V}}{\bigvee} \phi_{(v_i)}(p_1, ..., p_n).
\end{equation}

It is easy to see that for all sequences $(v_i)$, $\psi(v_1, ...v_n) = \phi_{(v_i)}(v_1, ..., v_n) = f(v_1, ..., v_n)$, since all of the $\phi$'s except for $\phi_{(v_i)}$ come out to $\bot$. It follows that, for all $(v_i)$, $\psi(v_1, ..., v_n) = f(v_1, ..., v_n)$. So $f$ is expressible as a propositional formula of $\lps$. \qed
\end{theorem}

Now we observe a number of nice rules which $\simp$ validates.

\begin{theorem}
The following are valid in $\lps$:
\begin{enumerate}
\item (Modus ponens.) $p \simp q, p \models q$.
\item (Contraction.) $\models (p \simp q) \simp ((p \simp (q \simp r)) \simp (p \simp r))$. 
\item (Weakened weakening.) $ \models p^\top \simp (q \simp p)$. 
\item (Contraposition.) $\models (p \simp q) \siff (\neg q \simp \neg p)$.
\item (Proof by cases.) $\models ((p^\top \simp q) \logand (p^+ \simp q) \logand (p^\bot \simp q)) \simp q$. 
\item (Proof by contradiction.) $\models (p \simp (q \logand q^\bot)) \simp p^\bot$. 
\end{enumerate}

\textbf{Proof.} These may be verifed using truth tables. \qed
\end{theorem}

$\lps$ faithfully preserves classical logic. This is true in at least two senses. Firstly, if an $\lps$ theory has only consistent models, then its $\lps$ consequences are precisely its classical consequences. This is trivial to verify, since in this case the $\lps$ models are precisely the classical models.

Secondly, if an $\lps$ theory $T$ has a consistent model, then we can add to it a non-trivializing axiom schema $U$ so that the $\lps$ consequences of $T \cup U$ are precisely the classical consequences of $T$. We let
\[ 
U = \{\forall x_1,...,x_n (R(x_1,...,x_n)^\con) : R\ \text{is an}\ n-\text{ary relation in the signature of}\ T\}.
 \]
$T \cup U$ has a model, since $T$ has a consistent model $M$ and $M \models U$. Clearly, furthermore, every model of $T \cup U$ is consistent; so the consequences of $T \cup U$ are just the classical consequences of $T \cup U$, which are just the classical consequences of $T$, since $U$ is tautologous in classical logic.\footnote{In translating $\lps$ formulas into classical logic, we translate $\simp$ as $\imp$. So $\phi^\con$ translates as $\phi \imp (\top \imp \phi)$.}

\section{Embeddings of classical theories}

An intended application of $\lps$ is in constructing inconsistent theories which prove all the theorems of some classical theory, and none of their negations. For instance, we might wish to construct an inconsistent set theory which nonetheless proves all theorems of ZFC and nothing false in ZFC. Our purpose in this section is to describe exactly what we mean by that. For this we introduce the notion of an ``embedding.''

Essentially, an embedding is a translation of the formulas of a classical theory $C$ into the formulas of an $\lps$ theory $L$. We map each relation $R$ in the signature of $C$ onto a corresponding predicate $\rho_R$ in the language of $L$. We let the domain of quantification of $C$ be translated as a subdomain of $L$'s domain of quantification characterized by a definable class $\kappa(x)$. This gives rise to a natural translation of the formulas of $C$ into the formulas of $L$. Then we require that $L$ proves all translated theorems of $C$, and additionally proves that all translated formulas are consistent. This implies that if $L$ is non-trivial, it does not prove the translations of any of the negations of the theorems of $C$.

\begin{definition}
Let $T$ be an $\lps$ theory, and $U$ a classical theory.\footnote{For simplicity we disallow constant symbols in the signature of $U$. They may be emulated using relation symbols.} An ``embedding'' of $U$ into $T$ consists of the following things:
\begin{enumerate}
\item For each $n$-ary relation symbol $R$ in the signature of $U$, an $\lps$ formula $\rho_R(x_1,...,x_n)$, with $x_1,...,x_n$ the free variables. This is the interpretation of the relation in the language of $U$.
\item A formula $\kappa(x)$ in one free variable $x$. This is intended to define a class of objects over which quantifiers range in the interpretation of $U$.
\end{enumerate}
Given a classical formula $\phi$, we define the $\lps$ translation $\tr(\phi)$ inductively as follows:
\begin{enumerate}
\item $\tr(R(x_1,...,x_n)) = \rho_R(x_1,...,x_n)$.
\item $\tr(\neg\phi) = \neg\tr(\phi)$.
\item $\tr(\phi \logand \psi) = \tr(\phi) \logand \tr(\psi)$.
\item $\tr(\forall x (\phi)) = \forall x (\kappa(x) \simp \tr(\phi))$.
\end{enumerate}
We require that for all formulas $\phi$ in the language of $U$, $T \models_\simp \tr(\phi)^\con$, and if $U \models_{\cl} \phi$ then $T \models_\simp \tr(\phi)$.
\end{definition}

\section{G\"odelian considerations}
\label{sec:godel}

In classical logic, consistent theories capable of expressing arithmetic form a partially ordered set where $T < U$ iff $U$ is capable of proving $T$ consistent. If $T < U$, then $T$ does not prove all the theorems of $U$. $\lps$ theories which embed classical theories factor into this hierarchy in a certain way: namely, a classical theory cannot prove nontrivival any \lps\ theory which embeds the classical theory.

In the following theorem we assume that an appropriate sound and complete deduction system has been defined for $\lps$ (which can readily be done). We write $\entails_\simp$ for $\lps$ syntactic entailment, and $\entails_\cl$ for classical syntactic entailment. Clearly if $T \entails_\simp \phi$ then there is an elementary arithmetical proof that $T \entails_\simp \phi$.

\begin{theorem}
\label{thm:nontriviality-hierarchy}
Let $T$ be an $\lps$ theory and $U$ a classical theory capable of expressing arithmetic, such that there is an embedding $E$ of $U$ into $T$, and $U$ proves that $E$ is an embedding. If $U$ is consistent, then $U$ does not prove $T$ nontrivial.

\textbf{Proof.} Suppose $U$ proves that $T$ is nontrivial. $U$ proves that if $U \entails_\cl \phi$ then $T \entails_\simp \tr(\phi)$. If $T \entails_\simp \phi \logand \phi^\bot$, then $T$ is trivial, and $U$ proves this. So $U$ proves that for all formulas $\phi$ in the language of $U$, $T \nvdash_\simp \tr(\phi) \logand \tr(\phi)^\bot$. So $U$ proves that for all classical $\phi$, $U \nvdash_\cl \phi \logand \phi^\bot$, i.e., $U \nvdash_\cl \phi \logand \neg\phi$. So $U$ proves its own consistency. \qed
\end{theorem}

\bibliography{ns-lp}

\end{document}